\documentclass{article}%
\usepackage{amsmath}
\usepackage{hyperref}
\usepackage{amsfonts}
\usepackage{amssymb}
\usepackage{graphicx}%
\providecommand{\U}[1]{\protect\rule{.1in}{.1in}}
\begin{document}

\title{\vspace{-0.5in}Elementary results for the fundamental representation of SU(3)}
\author{Thomas L. Curtright$^{\S }$ and Cosmas K. Zachos$^{\sharp}$\\Department of Physics, University of Miami\\Coral Gables, FL 33124-8046, USA}
\date{$^{\S }${\small curtright@miami.edu \ }$^{\sharp}${\small zachos@anl.gov}}
\maketitle

\begin{abstract}
\textit{A general group element for the fundamental representation of
}$SU\left(  3\right)  $\textit{\ can always be expressed as a second order
polynomial in a hermitian generating matrix }$H$,\textit{ with coefficients
consisting of elementary trigonometric functions dependent on the sole
invariant }$\det\left(  H\right)  $,\textit{ in addition to the group
parameter.}

\end{abstract}

\begin{center}
\textit{In memoriam Yoichiro Nambu (1921-2015)}
\end{center}

Consider an arbitrary $3\times3$ traceless hermitian matrix $H$. \ The
Cayley-Hamilton theorem \cite{Cayley} gives%
\begin{equation}
H^{3}=I\det\left(  H\right)  +\tfrac{1}{2}~H\operatorname*{tr}\left(
H^{2}\right)  \ , \label{C-H}%
\end{equation}
and therefore $\det\left(  H\right)  =\operatorname*{tr}\left(  H^{3}\right)
/3$. \ Note that an $H^{2}$ term is absent in the polynomial expansion of
$H^{3}$ because of the trace condition, $\operatorname*{tr}\left(  H\right)
=0$. \ Also note, since $\operatorname*{tr}\left(  H^{2}\right)  >0$ for any
nonzero hermitian $H$, this bilinear trace factor may be absorbed into the
normalization of $H$, thereby setting the scale of the group parameter space.

We may now write the exponential of $H$ as a matrix polynomial
\cite{Sylvester,Ilamed}. \ As a consequence of (\ref{C-H}) any such
exponential can be expressed as a matrix polynomial second-order in $H$, with
polynomial coefficients that depend on the displacement from the group origin
as a \textquotedblleft rotation angle\textquotedblright\ $\theta$. \ 

Moreover, the polynomial coefficients will also depend on invariants of the
matrix $H$. \ These invariants can be expressed in terms of the eigenvalues of
$H$, of course \cite{Sylvester,Ilamed,MSW,Rosen,Kusnezov,Laufer}.
\ Nevertheless, while the eigenvalues of $H$ will be manifest in the final
result given below, a deliberate diagonalization of $H$ is not necessary.
\ This is true for $SU\left(  3\right)  $, for a normalized $H$, since there
is effectively only one invariant: $\det\left(  H\right)  $. \ This invariant
may be encoded cyclometrically as another angle. \ Define%
\begin{equation}
\phi=\tfrac{1}{3}\left(  \arccos\left(  \tfrac{3}{2}\sqrt{3}\det\left(
H\right)  \right)  -\tfrac{\pi}{2}\right)  \ , \label{phi}%
\end{equation}
whose geometrical interpretation will soon be clear. \ Inversely,%
\begin{equation}
\det\left(  H\right)  =-\tfrac{2}{3\sqrt{3}}~\sin\left(  3\phi\right)  \ .
\label{det}%
\end{equation}
The result for any $SU\left(  3\right)  $ group element generated by a
traceless $3\times3$ hermitian matrix$\ H$ is then%
\begin{equation}
\exp\left(  i\theta H\right)  =\sum_{k=0,1,2}\left[  H^{2}+\tfrac{2}{\sqrt{3}%
}~H\sin\left(  \phi+\tfrac{2\pi k}{3}\right)  -\tfrac{1}{3}~I\left(
1+2\cos\left(  2\left(  \phi+\tfrac{2\pi k}{3}\right)  \right)  \right)
\right]  \frac{\exp\left(  \frac{2}{\sqrt{3}}~i\theta\sin\left(  \phi
+\frac{2\pi k}{3}\right)  \right)  }{1-2\cos\left(  2\left(  \phi+\frac{2\pi
k}{3}\right)  \right)  } \label{TheBigLebowski}%
\end{equation}
where we have set the scale for the $\theta$ parameter space by choosing the
normalization%
\begin{equation}
\operatorname*{tr}\left(  H^{2}\right)  =2\ . \label{normalization}%
\end{equation}
With this choice, the Cayley-Hamilton result (\ref{C-H}) is just \cite{Note}%
\begin{equation}
H^{3}=H+I\det\left(  H\right)  \ . \label{NormalizedC-H}%
\end{equation}
The normalization (\ref{normalization}) and the identity (\ref{NormalizedC-H})
are consistent with the Gell-Mann $\lambda$-matrices \cite{MGM}.\ 

So expressed as a matrix polynomial, the group element depends on the sole
invariant $\det\left(  H\right)  $ in addition to the group rotation
angle\ $\theta$. \ Both dependencies are in terms of elementary trigonometric
functions when $\det\left(  H\right)  $ is expressed as the angle $\phi$,
whose geometrical interpretation follows immediately from the three
eigenvalues of $H$ exhibited in the exponentials of (\ref{TheBigLebowski}).
\ Those eigenvalues are the projections onto three mutually perpendicular axes
of a single point on a circle formed by the intersection of the
$0=\operatorname*{tr}\left(  H\right)  $ eigenvalue plane with the
$2=\operatorname*{tr}\left(  H^{2}\right)  $ eigenvalue 2-sphere. $\ $The
angle $\phi$ parameterizes that circle. \ Equivalently, the eigenvalues are
the projections onto a single axis of three points equally spaced on a circle
\cite{Viete}.

Two cases deserve special mention. \ On the one hand, the Rodrigues formula
for $SO\left(  3\right)  $ rotations about an axis $\widehat{n}$, as generated
by $j=1$ spin matrices, is obtained for $\phi=0=\det\left(  H\right)  $.
\ Thus
\begin{equation}
\left.  \exp\left(  i\theta H\right)  \right\vert _{\phi=0}=I+iH\sin
\theta+H^{2}\left(  \cos\theta-1\right)  \ . \label{ER}%
\end{equation}
This is the Euler-Rodrigues result, upon identifying $H=\widehat{n}%
\cdot\overrightarrow{J}$\ (see \cite{CFZ,CvK}). \ It provides an explicit
embedding $SO\left(  3\right)  \subset SU\left(  3\right)  $. \ In fact,
(\ref{ER}) is true if $H$ is any \emph{one} of the first seven Gell-Mann
$\lambda$-matrices \cite{MGM}, or if $H$ is a normalized linear combination of
$\lambda_{1-3}$, or of $\lambda_{4-7}$. \ However, for generic linear
combinations of $\lambda_{1-7}$, $\det\left(  H\right)  $ will \emph{not}
necessarily vanish, and the general result (\ref{TheBigLebowski}) must be
used. \ 

On the other hand,
\begin{equation}
\lambda_{8}=\frac{1}{\sqrt{3}}\left(
\begin{array}
[c]{ccc}%
1 & 0 & 0\\
0 & 1 & 0\\
0 & 0 & -2
\end{array}
\right)
\end{equation}
is the only one among Gell-Mann's choices for the $3\times3$ representation
matrices for which $\phi\neq0$, and for which two eigenvalues are degenerate.
\ Obviously, $\det\left(  \lambda_{8}\right)  =\frac{-2}{3\sqrt{3}}$, so
$\phi=\pi/6$. \ In addition,%
\begin{equation}
\lambda_{8}^{2}=\frac{2}{3}~I-\frac{1}{\sqrt{3}}~\lambda_{8}\ .
\end{equation}
Thus, directly from (\ref{TheBigLebowski}),%
\begin{equation}
\exp\left(  i\theta\lambda_{8}\right)  =\tfrac{1}{3}\left(  2I+\sqrt{3}%
\lambda_{8}\right)  e^{\frac{1}{3}i\sqrt{3}\theta}+\tfrac{1}{3}\left(
I-\sqrt{3}\lambda_{8}\right)  e^{-i\frac{2}{\sqrt{3}}\theta}=\left(
\begin{array}
[c]{ccc}%
\exp\left(  i\theta/\sqrt{3}\right)  & 0 & 0\\
0 & \exp\left(  i\theta/\sqrt{3}\right)  & 0\\
0 & 0 & \exp\left(  -2i\theta/\sqrt{3}\right)
\end{array}
\right)  \ ,
\end{equation}
as it should. \ Note that this particular example followed from
(\ref{TheBigLebowski}) by \emph{carefully} taking the limit as $\phi
\rightarrow\pi/6$ of the $k=0$ and $k=1$ terms in that general expression (as
necessitated by the degeneracy of the corresponding eigenvalues of
$\lambda_{8}$) combined with the straightforward limit of the $k=2$ term.
\ That is to say,
\begin{gather}
\lim_{\phi\rightarrow\pi/6}\left(  \left[  \lambda_{8}^{2}+\tfrac{2}{\sqrt{3}%
}~\lambda_{8}\sin\left(  \phi+\tfrac{2\pi}{3}\right)  -\tfrac{1}{3}~I\left(
1+2\cos\left(  2\left(  \phi+\tfrac{2\pi}{3}\right)  \right)  \right)
\right]  \tfrac{\exp\left(  \frac{2}{\sqrt{3}}~i\theta\sin\left(  \phi
+\tfrac{2\pi}{3}\right)  \right)  }{1-2\cos\left(  2\left(  \phi+\frac{2\pi
}{3}\right)  \right)  }\right) \nonumber\\
=\lim_{\phi\rightarrow\pi/6}\left(  \left[  \lambda_{8}^{2}+\tfrac{2}{\sqrt
{3}}~\lambda_{8}\sin\left(  \phi\right)  -\tfrac{1}{3}~I\left(  1+2\cos\left(
2\phi\right)  \right)  \right]  \tfrac{\exp\left(  \frac{2}{\sqrt{3}}%
~i\theta\sin\phi\right)  }{1-2\cos\left(  2\phi\right)  }\right)  =\left(
\tfrac{1}{3}~I+\tfrac{1}{2\sqrt{3}}~\lambda_{8}\right)  e^{i\theta/\sqrt{3}%
}\ .
\end{gather}

Finally, one readily verifies that the Laplace transform of
(\ref{TheBigLebowski}) gives the resolvent in the standard form as a matrix
polynomial \cite{He,TSvK,TLC},%
\begin{equation}
\int_{0}^{\infty}e^{-t}\exp\left(  itsH\right)  dt=\frac{1}{I-isH}=\frac
{1}{\det\left(  I-isH\right)  }~\sum_{n=0}^{\operatorname*{rank}\left(
H\right)  -1}\left(  isH\right)  ^{n}\operatorname*{Trunc}%
_{\operatorname*{rank}\left(  H\right)  -1-n}\left[  \det\left(  I-isH\right)
\right]  \ , \label{Laplace}%
\end{equation}
where the truncation (as defined in \cite{CFZ}) is in powers of $s$. \ For the
case at hand $\operatorname*{rank}\left(  H\right)  =3$, again with
$\operatorname*{tr}\left(  H\right)  =0$ and $\operatorname*{tr}\left(
H^{2}\right)  =2$, so
\begin{equation}
\frac{1}{I-isH}=\frac{1}{1+s^{2}+is^{3}\det\left(  H\right)  }\left(  \left(
1+s^{2}\right)  I+isH-s^{2}H^{2}\right)  \ . \label{TheResolvent}%
\end{equation}
From this resolvent one immediately obtains a matrix polynomial for the simple
Cayley transform representation \cite{Cayley} of the corresponding $SU\left(
3\right)  $ group elements \cite{MSW},%
\begin{equation}
\frac{I+isH}{I-isH}=\frac{1}{1+s^{2}+is^{3}\det\left(  H\right)  }\left(
\left(  1+s^{2}-is^{3}\det\left(  H\right)  \right)  I+2isH-2s^{2}%
H^{2}\right)  \ . \label{CayleyTransform}%
\end{equation}
The Laplace transform (\ref{Laplace}) can be inverted, in standard fashion
\cite{Laplace}, to obtain (\ref{TheBigLebowski}) in terms of the impulse
response of the transfer function given by the prefactor in
(\ref{TheResolvent}). \ Explicitly,%
\begin{equation}
\exp\left(  i\theta H\right)  =\left(  H^{2}-iH\frac{d}{d\theta}-I\left(
1+\frac{d^{2}}{d\theta^{2}}\right)  \right)  \sum_{k=0,1,2}\frac{\exp\left(
\frac{2}{\sqrt{3}}~i\theta\sin\left(  \phi+2\pi k/3\right)  \right)  }%
{1-2\cos2\left(  \phi+2\pi k/3\right)  }\ .
\end{equation}

\paragraph{Acknowledgements}

This work was supported in part by NSF Award PHY-1214521, and in part by a
University of Miami Cooper Fellowship.


\begin{thebibliography}{99}                                                                                               %


\bibitem {Cayley}\textit{The Collected Mathematical Papers of Arthur Cayley},
Cambridge University Press (1889). \ See
\href{https://books.google.com/books?id=zJdQAAAAYAAJ&printsec=frontcover&source=gbs_ge_summary_r&cad=0#v=onepage&q&f=false}{Vol.
I.} pp 28-35 for the Cayley transform, and see
\href{https://books.google.com/books?id=TT1eAAAAcAAJ&printsec=frontcover&source=gbs_ge_summary_r&cad=0#v=onepage&q&f=false}{Vol.
II.} pp 475-496 for the Cayley-Hamilton theorem. \ Also see
\href{https://en.wikipedia.org/wiki/Cayley-Hamilton_theorem}{https://en.wikipedia.org/wiki/Cayley-Hamilton\_theorem}
and
\href{http://en.wikipedia.org/wiki/Cayley_transform}{https://en.wikipedia.org/wiki/Cayley\_transform}%
.

\bibitem {Sylvester}The general method to express \emph{any} analytic matrix
function of a finite, diagonalizable matrix as a polynomial in the matrix,
through the use of projection matrices, is due to J J Sylvester,
\href{http://babel.hathitrust.org/cgi/pt?id=mdp.39015024088141;view=1up;seq=287}{Phil.Mag.
16 (1883) 267-269}. \ Also see
\href{https://en.wikipedia.org/wiki/Sylvester's_formula}{https://en.wikipedia.org/wiki/Sylvester's\_formula}%
.

\bibitem {Ilamed}Y Lehrer, \textquotedblleft On functions of
matrices\textquotedblright%
\ \href{http://link.springer.com/article/10.1007/BF02848445}{Rendiconti del
Circolo Matematico di Palermo 6 (1957) 103-108};\newline Y Lehrer--Ilamed,
\textquotedblleft On the direct calculations of the representations of the
three-dimensional pure rotation group\textquotedblright%
\ \href{http://dx.doi.org/10.1017/S0305004100037452}{Proc.Camb.Phil.Soc. 60
(1964) 61--66} (especially see Remark (1), Eqn(10)).

\bibitem {MSW}A J MacFarlane, A Sudbery, and P H Weisz, \textquotedblleft On
Gell-Mann's $\lambda$-Matrices, $d$- and $f$-Tensors, Octets, and
Parametrizations of SU(3)\textquotedblright%
\ \href{http://link.springer.com/article/10.1007/BF01654302}{Commun.Math.Phys.
11 (1968) 77-90}.

\bibitem {Rosen}S P Rosen, \textquotedblleft Finite Transformations in Various
Representations of SU(3)\textquotedblright%
\ \href{http://scitation.aip.org/content/aip/journal/jmp/12/4/10.1063/1.1665634}{J.Math.Phys.
12 (1971) 673-681}.

\bibitem {Kusnezov}D Kusnezov, \textquotedblleft Exact matrix expansions for
group elements of SU(N)\textquotedblright%
\ \href{http://scitation.aip.org/content/aip/journal/jmp/36/2/10.1063/1.531165}{J.Math.Phys.
36 (1995) 898-906}.

\bibitem {Laufer}A Laufer, \textquotedblleft The exponential map of
GL(N)\textquotedblright%
\ \href{http://iopscience.iop.org/0305-4470/30/15/029/}{J.Phys.A:Math.Gen. 30
(1997) 5455}.

\bibitem {MGM}M Gell-Mann and Y Ne'eman,
\textit{\href{http://bookzz.org/book/1271076/4834ff}{\textit{The Eightfold
Way}, W A Benjamin (1964)}}. \ Also see
\href{https://en.wikipedia.org/wiki/Gell-Mann_matrices}{https://en.wikipedia.org/wiki/Gell-Mann\_matrices}%
.

\bibitem {Note}The reader may use (\ref{det}) and (\ref{NormalizedC-H}) to
check that the coefficients of the three exponentials in (\ref{TheBigLebowski}%
) are indeed projection matrices.

\bibitem {Viete}R W D Nickalls, \textquotedblleft Vi\`{e}te, Descartes and the
cubic equation\textquotedblright%
\ \href{http://www.jstor.org/stable/40378607}{Mathematical Gazette 90 (2006)
203--208}. \ Also see
\href{https://en.wikipedia.org/wiki/Cubic_function#Three_real_roots}{https://en.wikipedia.org/wiki/Cubic\_function\#Three\_real\_roots}%
.

\bibitem {CFZ}T L Curtright, D B Fairlie, and C K Zachos, \textquotedblleft A
Compact Formula for Rotations as Spin Matrix Polynomials\textquotedblright%
\ \href{http://www.emis.de/journals/SIGMA/2014/084/}{SIGMA 10 (2014) 084.}
\ e-Print: \href{http://arxiv.org/abs/1402.3541}{arXiv:1402.3541 [math-ph]}

\bibitem {CvK}T L Curtright and T S Van Kortryk, \textquotedblleft On
rotations as spin matrix polynomials\textquotedblright%
\ \ \href{http://iopscience.iop.org/1751-8121/48/2/025202/}{J.Phys.A:
Math.Theor. 48 (2015) 025202}. \ e-Print:
\href{http://arxiv.org/abs/1408.0767}{arXiv:1408.0767 [math-ph]}

\bibitem {He}M X He and P E Ricci, \textquotedblleft On Taylor's formula for
the resolvent of a complex matrix\textquotedblright%
\ \href{http://www.sciencedirect.com/science/article/pii/S0898122108003623}{Computers
and Mathematics with Applications 56 (2008) 2285--2288}.

\bibitem {TSvK}T S Van Kortryk, \textquotedblleft Cayley transforms of
$su\left(  2\right)  $ representations\textquotedblright\ e-Print:
\href{http://arxiv.org/abs/1506.00500}{arXiv:1506.00500 [math-ph]}

\bibitem {TLC}T L Curtright, \textquotedblleft More on Rotations as Spin
Matrix Polynomials\textquotedblright\ e-Print:
\href{http://arxiv.org/abs/1506.04648}{arXiv:1506.04648 [math-ph]}

\bibitem {Laplace}%
\href{https://en.wikipedia.org/wiki/Laplace_transform#Inverse_Laplace_transform}{https://en.wikipedia.org/wiki/Laplace\_transform\#Inverse\_Laplace\_transform}%

\end{thebibliography}
\end{document}